\documentclass[a4paper,12pt]{article}

\usepackage[T2A]{fontenc}
\usepackage[utf8]{inputenc}   % older versions of ucs package
\usepackage[russian]{babel}

\usepackage{CJKutf8}

\usepackage{amsmath,amstext,amsfonts,amsthm,amssymb}
\usepackage{eucal}
\usepackage{graphicx}
\renewcommand{\subsubsection}[1]{\addtocounter{subsubsection}{1}
{\ \\[3pt]\bf \thesubsubsection. \  #1} }
\swapnumbers

{  \theoremstyle{definition}

}
\newcommand{\Alt}{\operatorname{Alt}}
\newcommand{\an}{{\operatorname{an}}}

\newcommand{\Ch}{{\operatorname{Ch}}}

\newcommand{\Coker}{\operatorname{Coker}}
\newcommand{\Coul}{\operatorname{Coul}}

\newcommand{\GM}{{\operatorname{GM}}}

\newcommand{\Ker}{\operatorname{Ker}}
\newcommand{\KZ}{{\operatorname{KZ}}}

\newcommand{\Tot}{\operatorname{Tot}}

\newcommand{\hra}{\hookrightarrow}
\newcommand{\iso}{\overset{\sim}{\longrightarrow}}
\newcommand{\isom}{\overset{\sim}{=}}
\newcommand{\lra}{\longrightarrow}

\newcommand{\dpar}{\partial}

\newcommand{\bea}{\begin{eqnarray*}}
\newcommand{\eea}{\end{eqnarray*}}
\newcommand{\bean}{\begin{eqnarray}}
\newcommand{\eean}{\end{eqnarray}}

\newcommand{\teta}{\tilde\eta}

\newcommand{\tfg}{\tilde\fg}

\newcommand{\fb}{\mathfrak b}
\newcommand{\fg}{\mathfrak g}
\newcommand{\fh}{\mathfrak h}
\newcommand{\fm}{\mathfrak m}

\newcommand{\fn}{\mathfrak n}

\newcommand{\fsl}{\mathfrak{sl}}

\newcommand{\bm}{\mathbf{m}}
\newcommand{\bt}{\mathbf{t}}
\newcommand{\bz}{\mathbf{z}}

\newcommand{\CA}{\mathcal{A}}

\newcommand{\CC}{\mathcal{C}}

\newcommand{\CD}{\mathcal{D}}

\newcommand{\CH}{\mathcal{H}}
\newcommand{\CI}{\mathcal{I}}

\newcommand{\CL}{\mathcal{L}}
\newcommand{\CM}{\mathcal{M}}

\newcommand{\CO}{\mathcal{O}}
\newcommand{\CP}{\mathcal{P}}

\newcommand{\CZ}{\mathcal{Z}}

\newcommand{\BC}{\mathbb{C}}

\newcommand{\BN}{\mathbb{N}}

\newcommand{\BZ}{\mathbb{Z}}

\newcommand{\bze}{\text{\bf 0}}

\newcommand{\al}{\alpha}
\newcommand{\ka}{\kappa}
\newcommand{\om}{\omega}

\newcommand{\si}{\sigma}

\newcommand{\dl}{d\ln}

\newcommand{\ox}{\otimes}

\newcommand{\mc}{\mathcal}

\newcommand{\nc}{\newcommand}

\nc{\Id}{\text{Id}}
\nc{\la}{\lambda}

\begin{document}

\centerline{\bf DERIVED KZ EQUATIONS}

\bigskip\bigskip

\centerline{Vadim Schechtman\footnote{Institut de Math\'ematiques de Toulouse, Universit\'e Paul Sabatier, 118 route de Narbonne, 31062 Toulouse, France} and Alexander Varchenko\footnote{Department of Mathematics, University of North Carolina at Chapel Hill,
Chapel Hill, NC 27599-3250, USA; \ Faculty of Mathematics and
  Mechanics, Lomonosov Moscow State University,
Leninskiye Gory 1, 119991 Moscow GSP-1, Russia; \   Moscow Center of
Fundamental and Applied Mathematics
  Leninskiye Gory 1, 119991 Moscow GSP-1, Russia} }

\bigskip\bigskip

\begin{CJK}{UTF8}{min}

%\centerline{ヘッケの縮退代数とコ骸骨\footnote{hecke no shukutai daisu to kogaikotsu}}
%\centerline{ニルヘッケの代数とコ骸骨\footnote{niruhecke no daisu to kogaikotsu}}

\end{CJK}

%\bigskip\bigskip

 \centerline{December 27, 2020}

\

\

\hspace{9cm} {\it To the 30-th anniversary of} [SV]

\

\ 

\centerline{\bf Abstract}

\

In this note we strengthen the results of [SV] by presenting their derived version. 
Namely, we define a "derived Knizhnik - Zamolodchikov connection"\ and identify it with a "derived Gauss - Manin connection". 

\

\

\centerline{\bf \S 1. Introduction}

\bigskip

% REVIEW OF INTRODUCTION

\

{\bf 1.0. Brief review of the paper.} The main result of [SV] provided a realization of Knizhnik - Zamolodchikov equations arising in physics as equations on horizontal sections for a Gauss - Manin connection. 

More explicitly, without going into details to be given below, the KZ connection acts on a space 
of functions depending on $\bz\in B$ where $B$ is a domain in $\BC^n$ with values in a homology group $H_0(\fn, M)$ 
where $\fn$ is a certain Lie algebra, and $M$ a (maybe infinite dimensional) $\fn$-module. In other words the KZ connection 
acts on the trivial vector bundle over $B$ with a fiber $H_0(\fn, M)$, this vector bundle to be denoted 
$\CH_0(\fn, M)$. 

All homology spaces $H_i(\fn, M)$ are $\Lambda$-graded
$$
H_i(\fn, M) = \oplus_{\lambda\in \Lambda} H_i(\fn, M)_\lambda 
$$
where $\Lambda$ is certain lattice. For a given $\lambda$ only a finite number 
of spaces $H_i(\fn, M)_\lambda$, $0\leq i\leq N$, are different from $0$. Let us pick 
$\lambda$. 

On the other hand one has introduced in {\it op. cit.} a fibration (a smooth surjective map) 
$$
p_\lambda : X_\lambda \lra B
$$
and a $\CD$-module $\CL_\lambda$ over $X_\lambda$, and a finite group $\Sigma_\lambda$ 
(a product of symmetric groups) which acts on $X_\lambda$ and $\CL_\lambda$. 

One has constructed an isomorphism of the bundle  
$\CH_0(\fn, M)_\lambda$   equipped with 
the KZ connection with the bundle $(R^Np_{\lambda *}\CL_\lambda)^{\Sigma_\lambda}$ equipped with the GM connection.  

In fact in [SV] for all $0 \leq i\leq N$ there were established isomorphisms 
$$
\beta_{i,\lambda}:\ \CH_i(\fn, M)_\lambda \iso (R^{N-i}p_{\lambda *}\CL_\lambda)^{\Sigma_\lambda}
\eqno{(1.0.1)}
$$
of vector bundles over $B$. However the question of identification of a connection on 
$\CH_i(\fn, M)_\lambda$ corresponding to the GM connection on 
$(R^{N-i}p_{\lambda *}\CL_\lambda)^{\Sigma_\lambda}$ was left open for $i > 0$, although a natural 
candidate has been given. 

In the present note we establish this remaining point.  To do this we start from the remark 
that by its very definition in  
{\it op. cit.} isomorphisms (1.0.1) are induced by a map of complexes
$$
\eta_\lambda = (\eta_{\lambda,i}):\ \CC_\bullet(\fn, M)_\lambda \lra 
\Omega_{X_\lambda/B}^{N - \bullet}(\CL_\lambda)^{\Sigma_\lambda}
\eqno{(1.0.2)} 
$$
where $\CC_\bullet(\fn, M)_\lambda$ is the $\lambda$-homogeneous part of the Chevalley chain complex, and  $\Omega_{X_\lambda/B}^{\bullet}(\CL_\lambda)$ is 
certain complex of differential form on $X_\lambda$, the relative de Rham complex of 
$\CL_\lambda$.   

A naive expectation would be that:

$(a)$ {\it for the KZ part}: 

each term $\CC_i(\fn, M)_\lambda$ comes equipped with an  integrable connection, these 
connections are compatible with differentials and thus induce a connection on the 
cohomology $\CH_i(\fn, M)_\lambda$;   

$(b)$ {\it for the GM part}:   

similarly, each term $\Omega_{X_\lambda/B}^{j}(\CL_\lambda)$ comes equipped with an 
integrable $\Sigma_\lambda$-equivariant connection, these 
connections are compatible with differentials and thus induce a connection on the 
cohomology  $R^{j}p_{\lambda *}\CL_\lambda$; 

$(c)$  the map $\eta_\lambda$ is compatible with the connections in (a), (b), and therefore the  
isomorphisms $\beta_{i,\lambda}$  (1.0.1) identify two connections.

In reality, $(a)$ is literally true (and easy); this is present in [SV].

Point $(b)$ is more delicate: there is {\it no} natural connection on the complex 
$\Omega_{X_\lambda/B}^{\bullet}(\CL_\lambda)$. Happily, to define a connection 
on the cohomology a weaker structure is sufficient: 

$(b')$ there exists a filtered complex  
such that the term $E_1$ of the corresponding spectral sequence (recalled in Appendix) coincides with the de Rham 
complex of the GM connection on  $R^{\bullet}p_{\lambda *}\CL_\lambda$. 

This filtered complex is described below: it is a generalization of the Katz - Oda construction for the GM connection, [KO]. 

Accordingly, $(c)$ should be replaced by

$(c')$ the map $\eta_\lambda$ may be extended to a map of filtered complexes which, after 
passing to $E_1$-terms, induces a map from the de Rham complex of the KZ connection 
to the de Rham complex of the GM connection.

\ 

Now we will describe some details of what was said above.

{\bf 1.1. Knizhnik - Zamolodchikov connection.} Let $\fg$ be a  complex 
Lie algebra equipped with an element 
$$
\Omega \in \fg\otimes\fg
$$
having the following property:

{\bf 1.1.1.} {\it Let $M_1, M_2$ be arbitrary $\fg$-modules. 
The actions of $\Omega$ and $\fg$ on $M_1\otimes M_2$ commute.}

{\bf 1.1.2. Example.}  Let $\fg$ be finite dimensional, equipped with a non-degenerate invariant symmetric bilinear form
$( , )$. Denote 
$$
\Omega = \sum_i x_i\otimes x^i\in \fg\otimes\fg
$$
where $\{x_i\}\subset \fg$ is any $\BC$-base, and $\{x^i\}$ is the dual base, i.e. 
$(x_i, x^i) = \delta_{ij}$. This element ("the Casimir") does not depend on a choice of a base and satisfies 1.1.1. 

% {\bf Proof.} Adapt [S], Part I, Ch. VI, \S 3, proof of H.Weyl's theorem. (OR BOURBAKI?)

Let $M_1, \ldots , M_n$ be  $\fg$-modules, $n\geq 1$. Denote 
$M = M_1\otimes \ldots\otimes  M_n$. 

\

For a smooth affine complex\footnote{in what follows the base field $\BC$ of complex numbers may be replaced by any  field of characteristics $0$} algebraic variety $U$, $\Omega^\bullet(U)$ will denote the space of global sections for its   
{\it algebraic} de Rham complex $\Omega^\bullet_U$. Thus $\Omega^0_U = \CO_U$ is the sheaf of functions, etc.

If $M$ is a vector space, we denote 
$$
\Omega^\bullet(U; M) := \Omega^\bullet(U)\otimes M.
$$

Let $n\geq 1$ be an integer. Let $M_1, \ldots , M_n$ be $\fg$-modules; set 
$M = M_1\otimes \ldots\otimes  M_n$. For each $i\neq j$ we have an operator
$$
\Omega_{ij}: M \lra M
$$
acting as $\Omega$ on $M_i\otimes M_j$ and as identity on the other factors. 

Denote
$$
U_n = \{\bz = (z_1, \ldots, z_n) \in \BC^n|\   z_i\neq z_j\ \text{for all}\ i\neq j\}
\eqno{(1.1.1)}
$$
Thus $U_1 = \BC$. 

The  KZ connection is an operator 
$$
\nabla_{KZ}: \Omega^0(U_n; M) = \CO(U_n)\otimes M\lra \Omega^1(U_n; M)
$$
given by 
$$
\nabla_{KZ} = d_{DR} + \Omega_{KZ} := d_{DR} - \frac{1}{\kappa}\sum_{i < j}\frac{\Omega_{ij}(dz_i - dz_j)}{z_i - z_j}
\eqno{(1.1.2)}
$$
where $d_{DR}$ is the de Rham differential. Here $\kappa\in \BC^*$ is a complex parameter. 

Thus $\nabla_{KZ} = d_{DR}$ if $n = 1$. 

This connection is integrable: if we define, starting from $\nabla_{KZ}$, operators
$$
\nabla_{KZ}:\ \Omega^i(U_n; M)\lra \Omega^{i+1}(U_n; M)
$$
for all $i$ in the usual way then $\nabla_{KZ}^2 = 0$ (this amounts to the classical YB equation for the differential form $\Omega_{KZ}$).

In other words, $\nabla_{KZ}$ is an integrable connection (i.e. it defines a structure of a 
$\CD_{U_n}$-module) on  the trivial bundle $\CM$ over $U_n$ 
with fiber $M$.

{\bf 1.2. The Chevalley complex and the derived KZ.} Let $\fn\subset\fg$ be a Lie subalgebra.

We will be interested in Chevalley chain complexes
$$
C_\bullet(\fn, M):\ \ldots \lra \Lambda^2\fn\otimes M \lra \fn\otimes M \lra M \lra 0
$$
where $d(g\otimes x) = gx$, 
$$
d(g_1\wedge g_2\otimes x) = g_1\otimes g_2x - g_2\otimes g_1x - [g_1,g_2]\otimes x,
$$
etc.

Let $\CC_\bullet(\fn, M)$ denote the trivial vector bundle over $U_n$ with a fiber 
$C_\bullet(\fn, M)$, so it is a complex of vector bundles.

\

We define the {\it derived KZ connection} as an integrable connection on 
$\CC_\bullet(\fn, M)$ given by the same formula as above,
$$
\nabla_{KZ} = d_{DR} + \Omega_{KZ} := d_{DR} - \frac{1}{\kappa}\sum_{i < j}\frac{\Omega_{ij}(dz_i - dz_j)}{z_i - z_j}
\eqno{(1.2.1)}
$$
where now the operators
$$
\Omega_{ij}:\ C_l(\fn, M) = \Lambda^l\fn\otimes M \lra C_l(\fn, M)
$$
are acting through the factor $M$. 

Whence we get the corresponding de Rham complex
$$
\Omega_{KZ}^\bullet(U_n, C_\bullet(\fn, M)) = DR(\CC_\bullet(\fn, M), \nabla_{KZ})(U_n). 
\eqno{(1.2.2)}
$$
It is a  double complex: the commutation of the Chevalley differential with 
$\nabla_{KZ}$ follows from 1.1.1. 

We call it the {\it KZ-Chevalley complex}. 

%DEGREES

% ACTUALLY THEY APPEAR ALREADY IN 

In fact this complex appears {\it avant la lettre} already in [SV] 7.2.3.

{\bf 1.3. Derived Gauss - Manin connection.} Let $N\geq 0$ be an integer. Consider the affine space $\BC^{n+N}$ with coordinates $z_1, \ldots, z_n, t_1, \ldots, t_N$, and inside it an open subspace
$$
U_{n,N} = \{z_i \neq z_j,\ z_i \neq t_a,\ t_a\neq t_b\}.
$$
We have an obvious projection
$$
p: U_{n,N} \lra U_n.
$$
The de Rham algebra $\Omega^\bullet(U_{n,N})$ is the total complex of 
a bicomplex
$$
\Omega^\bullet(U_{n,N}) = \Tot\Omega^{\bullet\bullet}(U_{n,N})
$$
where $\Omega^{pq}(U_{n,N})$ is the space of forms containing $p$ differentials $dt_i$ and 
$q$ differentials $dz_m$, the full de Rham differential being the sum
$$
d_{DR} = d_z + d_t.
$$
The {\it relative} de Rham complex is by definition
$$
\Omega^\bullet(U_{n,N}/U_n) = (\Omega^{0\bullet}(U_{n,N}), d_t); 
$$
one has a projection
$$
p:\  \Omega^{\bullet}(U_{n,N}) \lra \Omega^\bullet(U_{n,N}/U_n)
$$
Let $\CL$ be a $\CD_{U_{n,N}}$-module, i.e. a quasicoherent $\CO_{U_{n,N}}$-module equipped 
with an integrable connection
$$
\nabla:\ \CL \lra \Omega^1_{U_{n,N}}\otimes \CL;
$$
its de Rham complex is
$$
DR(\CL):\ 0 \lra  \CL \overset{\nabla}\lra \Omega^1_{U_{n,N}}\otimes \CL \overset{\nabla}\lra \Omega^2_{U_{n,N}}\otimes \CL \lra \ldots
$$ 
By definition, {\it the de Rham complex of the derived Gauss - Manin connection} on the direct image $Rp_*\CL$, to be denoted $DR(Rp_*\CL)$, is the same complex $DR(\CL)$ equipped with a decreasing filtration      
$$
F_z^0DR(\CL) = DR(\CL)\supset F_z^1DR(\CL)\supset \ldots
\eqno{(1.3.1)}
$$
where $F_z^i DR(\CL)$
is the subcomplex containing $\geq i$ differentials $dz_a$.

Note that the utmost left column of $F^0/F^1$ is the 
relative de Rham complex representing $Rp_*\CL$ whose cohomology are 
the sheaves $R^ip_*\CL$. These sheaves carry the usual GM connections $\nabla^i$.

The complexes $E^i(DR(Rp_*\CL), F^\bullet_z)$ defined in the Appendix, A2.1 (the components of the $E_1$ term 
of the spectral sequence for our filtered complex) are nothing else but their de Rham complexes of $R^ip_*\CL$: 
$$
E^i(DR(Rp_*\CL), F^\bullet_z) \isom DR(R^ip_*\CL, \nabla^i).
$$
This isomorphism justifies the above definition. 

{\bf 1.3.1.  Remark.} For the case of a trivial connection on $\CO_{U_{n,N}}$ the above construction is nothing else but 
the Katz - Oda definition of the usual GM connection, [KO]. 

\

{\bf 1.4. Coulomb $\CD$-modules.} 
Let $V$ be a finite dimensional complex vector space equipped with a symmetric bilinear form 
$( , )$. Let 
$$
\mu = (\mu_1, \ldots, \mu_n)\in V^n,\ 
\alpha = (\alpha_1, \ldots, \alpha_N)\in V^N,
$$
and $\kappa\in \BC^*$.

We associate to these data a $\CD$-module $\CL(\mu, \alpha)$, to be called {\it a Coulomb
\footnote{"Loi fondamentale de l'\'Elictricit\'e. La force r\'epulsive des deux petits globes \'electris\'es de la même nature d'électricit\'e, est en raison inverse du carré de la distance du centre de deux globes."\ Charles-Augustin de Coulomb, {\it Premier M\'emoire sur l’\'Electricit\'e et le Magn\'etisme}, 1785.} $\CD$-module}, over $U_{n,N}$: by definition 
it is the structure sheaf $\CO_{U_{n,N}}$ equipped with a connection
$$
\nabla(\mu, \alpha) = d_{DR} + \frac{1}{\kappa}\omega(\mu, \alpha)
$$
where 
$$
\omega(\mu, \alpha) = \sum_{i < j} (\mu_i, \mu_j)d\ln(z_i - z_j) - 
$$
$$
- \sum_{i, a} (\mu_i, \alpha_a)d\ln(z_i - t_a) + 
\sum_{a < b} (\alpha_a, \alpha_b)d\ln(t_a - t_b).
$$

{\bf 1.5.} On the other hand we can associate with the data $(V, \mu, \alpha)$ above 
a Lie algebra $\fg = \fg(\alpha)$ ("a Kac-Moody algebra without Serre relations") and a collection of "contragradient Verma"\ $\fg$-modules 
$M(\mu_1)^c, \ldots, M(\mu_n)^c$. 

For example if $V$ is one-dimensional and $\alpha_1 = \ldots = \alpha_N$ then 
$\fg = \fsl_2$. 

Let 
$$
M = M(\mu_1)^c\otimes \ldots \otimes  M(\mu_n)^c,
$$
and consider the total complex of the de Rham complex (1.2.2) 
$\Tot\Omega^\bullet(U_n,C_\bullet(\fn, M))$. It is $\Lambda$-graded where 
$\Lambda = \sum_i \BZ\alpha_i\subset V$, and it carries a decreasing filtration 
$$
F^\bullet_z\Tot\Omega^\bullet(U_n,C_\bullet(\fn, M))
$$
where
$$
F^\bullet_z\Tot\Omega^\bullet(U_n,C_\bullet(\fn, M))\subset 
\Tot\Omega^\bullet(U_n,C_\bullet(\fn, M))
$$
is the subcomplex of differential forms containing $\geq i$ differentials $dz_a$.

Let $\lambda = \sum_i \alpha_i\in \Lambda$.
Our main result defines a map from the  $\lambda$-homogeneous component of this filtered complex to the filtered complex $(DR(Rp_*\CL(\mu, \alpha), F^\bullet_z)$.

For details see Theorem 3.8 and Corollary 3.9.

\

{\it Plan of the paper}

\

In the next $\S 2$ we discuss in detail the case $\fg = \fsl_2$. The general case 
is discussed in $\S 3$. In the Appendix we recall some standard homological algebra of filtered complexes. 

\

{\bf 1.6. Acknowledgements.} We are grateful to B.Toen and D.Gaitsgory for useful conversations.  
A. Varchenko was supported in part by NSF grant DMS-1954266.

\

\

\

\

\centerline{\bf \S 2. The case $\fg = \fsl_2$}

\

{\bf 2.0. Setup.} We consider the Lie algebra $\fg = \fsl_2$ with standard generators 
$e, f, h$; let $\fn := \BC f\subset \fg$ (resp. $\fn_+ :=\BC e$) be the lower (resp. upper) triangular subalgebra. We will identify $\fn_+$ with $\fn^*$, with $e$ being dual to $f$. 

The Casimir element is
$$
\Omega = \frac{1}{2}h\otimes h + e\otimes f + f\otimes e. 
$$

{\bf 2.0.1. Invariance lemma.} {\it Let $M_1, M_2$ be arbitrary $\fg$-modules. 
The actions of $\Omega$ and $\fg$ on $M_1\otimes M_2$ commute.}

{\bf Proof}: exercise for the reader.

% NOT FORGET TO MENTION THIS FOR AN ARBITRARY $\fg$

\

{\bf 2.1. Chevalley complex.}

If $M$ is a $\fg$-module, $C_\bullet(\fn^*, M^*) = C_\bullet(\fn_+, M^*)$ will denote the Chevalley chain complex
$$
0 \lra \fn^*\otimes M^* \overset{d^*}\lra M^* \lra 0
\eqno{(2.1.1)}
$$
living in degrees $-1, 0$. Here the action of $\fn^*$ on the dual space $M^*$ is given by 
$$
(f^*\alpha)(x) = \alpha(ex),\ x\in M,\ \alpha\in M^*
\eqno{(2.1.1a)}
$$
where $f^*\in \fn^*$ is defined by $f^*(f) = 1$. 
 
Next, 
$C^\bullet(\fn, M^c)$ will denote the dual complex
$$
0\lra M \overset{d}\lra \fn\otimes M \lra 0,\ d(x) = f\otimes ex
\eqno{(2.1.2)}
$$
living in degrees $0,1$. 

For $m\in \BC$, $M(m)$ will denote the Verma module with a vacuum vector $v$ such that 
$hv = mv, ev = 0$. It is $\BN$-graded:
$$
M(m) = \oplus_{k\geq 0} M(m)_k
$$
where $M(m)_k = \BC f^kv$.

Fix a natural $n\geq 1$ and an $n$-tuple $\bm = (m_1, \ldots, m_n)\in \BC^n$, and consider the tensor product
$$
M(\bm) = M(m_1)\otimes\ldots\otimes M(m_n)
$$
The above grading on each $M(m_i)$ gives rise to an $\BN$-grading on $M(\bm)$:
$$
M(\bm) = \oplus_{k\geq 0} M(\bm)_k
$$
where
$$
M(\bm)_k = \oplus_{k_1+\ldots+k_n = k} M(m_1)_{k_1}\otimes\ldots\otimes M(m_n)_{k_n}.
$$
For a multi-index $a = (k_1, \ldots, k_n)$ we denote
$$
|a| = \sum_{i=1}^n\ a_i,
\eqno{(2.1.3)}
$$
and
$$
f^av := f^{a_1}v_1\otimes\ldots\otimes f^{a_n}v_n\in M(\bm)_{|a|}.
\eqno{(2.1.4)}
$$
The Chevalley complex acquires a grading as well: 
$$ 
C^\bullet(\fn, M(\bm)^c) = \oplus_{k\geq 0} C^\bullet(\fn, M(\bm)^c)_k,
$$
with
$$
C^\bullet(\fn, M(\bm)^c)_k: 0 \lra  M(\bm)_k \lra \fn\otimes M(\bm)_{k-1}\lra 0.
\eqno{(2.1.5)}  
$$

% COROLLARY: DERIVED KZ IS CORRECTLY DEFINED.

\

{\bf 2.2. Logarithmic forms.} Recall that $\kappa\in \BC^*$ is fixed.

We fix an integer $N\geq 0$ and consider the space 
$U_{n,N}$ (see 1.4 above).

We are going to define certain logarithmic forms on this space. For a function $u$ we denote
$$
d\ln u := \frac{du}{u}
$$
The symmetric group $\Sigma_N$ acts on  forms from $\Omega^\bullet(U_{n,N})$ by permuting variables $ t_1,\dots, t_N$.

For a  differential form $w$ we define by $\Alt \,w$ the skew-symmetrization of $w$ with respect to the $\Sigma_k$-action,
\bea
\Alt\, w(t_1,\dots,t_N) = \sum_{\si\in \Sigma_N} (-1)^{\si}w(t_{\si(1)},\dots,t_{\si(N)}).
\eea

All forms appearing in our constructions are skew-symmetric. They are given by the following formulas.
For $a=(a_1,\dots,a_n)\in \BN^n$, $|a| := \sum a_i =N$, we define 
\bea
w_a =\frac 1{a_1!\dots a_n!}\Alt\, u_a,
\eea
 where
\bea
u_a
&=& \dl(t_1-z_1)\wedge\dots\wedge  \dl(t_{a_1}-z_1) + 
\\
&&
+ \,
\dl(t_{a_1+1}-z_2)\wedge\dots \wedge \dl(t_{a_1+a_2}-z_2) +  
\\
&&
\dots
+\,
\dl(t_{a_1+\dots+a_{n-1}+1}-z_n)\wedge\dots \wedge \dl(t_{k}-z_n).
\eea

Similarly, for $b=(b_1,\dots,b_n)$, $|b|= N - 1$, we define 
\bea
w_b =\frac 1{b_1!\dots b_n!}\Alt\, u_b,
\eea
 where
\bea
u_b
&=& - \kappa \biggl( \,\
\dl(t_2-z_1)\wedge\dots\wedge  \dl(t_{b_1+1}-z_1) + 
\\
&&
+\,
\dl(t_{b_1+2}-z_2)\wedge\dots \wedge \dl(t_{b_1+b_2+1}-z_2) +  
\\
&&
\dots
+ \,
\dl(t_{b_1+\dots+b_{n-1}+2}-z_n)\wedge\dots \wedge \dl(t_{N}-z_n)\biggr).
\eea

In this formula we start from the variable $t_2$ and have the factor $ -\ka$ in front of the exterior product.

For example if $N=2$, $a=(2,0)$, $b=(1,0)$, then
\bea
w_a
&=&
 \dl(t_1-z_1)\wedge \dl(t_2-z_1)  
\\
w_b
&=&
- \kappa(\dl(t_2-z_1) + \dl(t_1-z_2)).
\eea

\

{\bf 2.3. Coulomb $\CD$-module.} Define a "Coulomb interaction"\ closed  $1$-form
$$
\om_{\bm} :=
\sum_{1\leq s<u\leq n}\frac{m_sm_u}2\,\dl(z_s-z_u)+
\sum_{1\leq i<j\leq N}2\,\dl(t_i-t_j) - 
$$
$$
- \sum_{i=1}^N\sum_{s=1}^n m_s\,\dl(t_i-z_s)\in \Omega^1(U_{n,N})
\eqno{(2.3.1)}
$$
Define a differential $\nabla_\bm$ on the graded space $\Omega^\bullet(U_{n,N})$
$$
\nabla_\bm := d_{DR} + \frac{1}{\kappa}\omega_\bm: \  \Omega^i(U_{n,N}) \lra \Omega^{i+1}(U_{n,N}) 
$$
Note that $\nabla_\bm^2 = 0$ since $d_{DR}\omega_\bm = 0$.

We will denote by $\Omega_\bm^\bullet(U_{n,N})$ the space $\Omega^\bullet(U_{n,N})$ equipped with the differential $\nabla_\bm$. 

This is nothing else but the complex of global sections 
 for  the de Rham complex $DR(\CL(\bm,N))$ of the 
{\it Coulomb $\CD$-module $\CL(\bm) = \CL(\bm,N)$} over $U_{n,N}$ which is by definition 
the structure sheaf 
$\CO_{U_{n,N}}$ equipped with a connection $\nabla_\bm := d_{DR} + \frac{1}{\kappa}\omega_\bm$.

\ 

{\bf 2.4. Coulomb - KZ - Chevalley complex and a canonical $N$-cocycle in it.} Recall a Chevalley complex $C^\bullet(\fn, M(\bm)^c)_N$.

Consider a double complex which as a bigraded vector space is a tensor product
$$
C^{\bullet\bullet}_{\bm, N} := \{ C^{pq}_{\bm, N}\}  
$$
where
$$
C^{pq}_{\bm, N} := \Omega^p(U_{n,N}) \otimes C^q(\fn, M(\bm)^c)_N 
$$
Note that along $q$-axis it has only two nontrivial components: $0 \leq q\leq 1$. 

By definition it is equipped with two differentials: 

--- the horizontal one is a KZ - Coulomb differential
$$
\nabla_{\KZ, \Coul} = d_{DR} + \frac{1}{\kappa}\omega_\bm - \frac{1}{\kappa}\omega_{\KZ}  
$$
where
$$
\omega_{\KZ} := \sum_{1\leq i < j \leq n}\Omega_{ij}d\ln(z_i - z_j)
\eqno{(2.4.1)}
$$
% SIGNS ?

It acts on the index $p$:
$$
\nabla_{\KZ, \Coul}: \Omega^p(U_{n,N}) \otimes C^q(\fn, M(\bm)^c)_N \lra 
\Omega^{p+1}(U_{n,N}) \otimes C^q(\fn, M(\bm)^c)_N
$$

--- the vertical one is the Chevalley differential $d_{\Ch}$ acting on the second factor. 

We will be interested in the associated total complex
$$
C^{\bullet}_{\bm, N} := \Tot C^{\bullet\bullet}_{\bm, N}.
$$
Recall the notations (2.1.4).

Define elements
\bea
\CI_0 
&:=&
 \sum_{|a|=N} w_a\otimes f^{(a)}v\in C^{N0}_{\bm, N}
\\
\CI_1
&:=&
 \sum_{|b|=N-1} w_b\otimes (f\otimes f^{(b)}v)\in C^{N-1,1}_{\bm, N}
 \\
 \CI
 &:=&
 \CI_0 + \CI_1\in  C^{N}_{\bm, N}
 \eea

%SIGNS

{\bf 2.5. Theorem.} {\it $\CI$ is a cocycle in $C^{\bullet}_{\bm, N}$ of total degree $N$. In components:
$$
\nabla_{\KZ, \Coul}\CI_0 = 0
\eqno{(2.5.1)}
$$
$$
d_{\Ch}\CI_0 + \nabla_{\KZ, \Coul}\CI_1 = 0.
\eqno{(2.5.2)}
$$}

{\bf Proof.}
We deduce Theorem 2.5 from the two main results in [SV].  The first of them is 
[SV, Theorem 6.16.2] on the relation between the Lie algebra differential and the de Rham differential.
 The second  is [SV, Theorem $7.2.5''$] on the relation between the KZ equations and the Gauss-Manin connection.
 
Since all forms $w_a$ are closed 
 the equation (2.5.1) may be rewritten as
 $$
 \frac{1}{\kappa}\biggl(\omega_\bm - \omega_{KZ}\biggr)\CI_0 = 0.
 $$ 
 This equation is the statement of [SV, Theorem 7.2.5''] applied to the $\frak{sl}_2$ case.

 Equation  (2.5.2) may be rewritten as
 $$
 \frac{1}{\kappa}\biggl(\omega_\bm - \omega_{KZ}\biggr)\CI_1 + d_{Ch}\CI_0 = 0
 $$
 and it  can be split into two equations.
 
 One of these equations  follows      from [SV, Theorem 7.2.5'']
 applied to the situation with $N-1$ of $t$-variables instead of the $N$ variables $t_1,\dots, t_N$,
 and the other equation follows from [SV, Theorem 6.16.2].

More precisely, consider the splitting
$$
% \bean
% \label{3}
 \frac{1}{\kappa}\omega_\bm
  \CI_1 =\CP_1+\CP_2,
% \eean
 $$
 where $\CP_1, \CP_2$ are defined as follows.   
 We have 
 \bea
 \CI_1 = \sum_{|b|=N-1}\sum_{\si\in S_N} (-1)^\si u_b(t_{\si(2)},\dots,t_{\si(N)}) 
\ox\Big(f\otimes f^{(b)}v\Big)
\eea 
  and $\om_\bm$ is  the sum of 1-forms, $\om_\bm =\sum_\al \om_\al$, see (2.3.1).
   We say that a summand
 \bea
 (-1)^\si \om_\al   \wedge  u_b(t_{\si(2)},\dots,t_{\si(k)}) \ox\Big(f\otimes f^{(b)}v\Big)
  \eea 
  belongs to $\mc P_1$ if $\om_\al$ does not have the variable $t_{\si(1)}$, otherwise it belongs to $\mc P_2$.

{\bf 2.5.1. Lemma.} {\it  We have
$$
\mc P_1 -\frac 1\ka \omega_{\KZ}\CI_1=0,
\eqno{(2.5.1.1)}
$$
$$
\mc P_2 +d_{\Ch}\CI_0 = 0.
\eqno{(2.5.1.2)}
$$}
  
{\bf Proof} of the Lemma. 
 Equation (2.5.1.2) follows from [SV], Theorem 6.16.2. Equation
(2.5.1.1) follows from [SV], Theorem $7.2.5''$.
$\square$

This implies (2.5.2) and  achieves the proof of 2.5. $\square$

\

{\bf 2.6. Interpretation of the cocycle $\CI$ as a map $\eta: DR($KZ) $\lra DR$(GM).} 
  Note that the Coulomb de Rham complex $\Omega^\bullet_\bm(U_{n,N})$  is a dg-module over the de Rham algebra $\Omega^\bullet(U_{n,N})$ which in turn is    
  a $\Omega^\bullet(U_{n})$-algebra due to the projection $p: U_{n,N}\lra U_n$. 
  
  %(this is an instance of a general fact that $\CD$-modules are\  
  %"the same as"\ dg-modules over the de Rham algebra, cf. [K]\footnote{an example of Koszul duality}), 
  
 Consider the trivial vector bundle $\CM(\bm)$ over $U_n$ with a fiber $M(\bm)$; it carries the integrable KZ connection 
  $$
  \nabla_{\KZ} = d_z - \frac{1}{\kappa}\omega_{\KZ}
  \eqno{(2.6.1)}
  $$
  which makes of it a $\CD_{U_n}$-module. The space of global sections of its de Rham complex will be
  $$
   DR(\CM(\bm))(U_n) = \Omega^\bullet(U_n)\otimes_\BC M(\bm).
  $$
As usual this object is $\BN$-graded.

Next we can pass to Chevalley chains and consider a complex of vector bundles 
$$
\CC_\bullet(\fn^*, M(\bm)^*_N) = C_\bullet(\fn^*, \CM(\bm)^*_N),  
$$
whose dual will be
$$
\CC^\bullet(\fn, M(\bm)^c_N) = C^\bullet(\fn, \CM(\bm)^c_N).  
$$
Both complexes carry  KZ connections induced by (2.7.1); therefore we may consider their 
de Rham complexes which are $\Omega^\bullet_{U_n}$-modules. 

% CHECK $\fn$ VERSUS $\fn^*$
  
  Our main hero, the KZ - Coulomb - Chevalley complex may be rewritten in a form
$$
C^{\bullet\bullet}_{\bm, N} = DR(\CC^\bullet(\fn, M(\bm)^c_N)(U_n)
\otimes_{\Omega^\bullet(U_n)}
\Omega^\bullet_\bm(U_{n,N})    
$$
By linear algebra, to give a $0$-cocycle 
  $$
  Z\in \Tot(A^\bullet\otimes B^\bullet)^0
  $$ 
  in the total complex of a tensor product of two complexes $A^\bullet\otimes B^\bullet$ is equivalent to giving a map of complexes 
  $$
  \eta(Z):\ A^{\bullet *}\lra B^\bullet.
  $$ 
  
  Therefore our cocycle $\CI$ gives rise to a map between two complexes
  $$
  \eta = \eta(\CI): DR(\CC_\bullet(\fn^*, M(\bm)^*_N)(U_n) \lra \Omega^\bullet_\bm(U_{n,N})[N] 
\eqno{(2.6.2)}
  $$
Both complexes are filtered: 

namely, we define 
$$
F^iDR(\CC_\bullet(\fn^*, M(\bm)^*_N)(U_n)\subset DR(\CC_\bullet(\fn^*, M(\bm)^*_N)(U_n)
$$
to be the subcomplex of forms of degree $\geq i$, and
$$
F_z^i\Omega^\bullet_\bm(U_{n,N}) \subset \Omega^\bullet_\bm(U_{n,N})
$$
to be the subcomplex of forms containing $\geq i$ differentials $dz_a$. 

{\bf 2.6.1. Key fact.} {\it The map $\eta$ is compatible with the filtrations.} $\square$

As a corollary , the induced map of $E_1$-terms of the corresponding spectral sequences gives rise to maps between the de Rham complexes 
  $$
  \eta^i:\ DR(H_i(\fn^*, \CM(\bm)^*_N), \nabla_\KZ) \lra 
  DR(Rp_*^{N-i}\CL(\bm, N), \nabla_\GM),
  $$
  $0\leq i\leq 1$, cf 1.3. 
  
% SYMMETRIC INVARIANTS
 
 By construction these maps land in the subsheaves of anti-invariants
$$
  \eta^i:\ DR(H_i(\fn^*, \CM(\bm)^*_N), \nabla_\KZ) \lra 
  DR(Rp_*^{N-i}\CL(\bm, N)^{\Sigma_N,-}, \nabla_\GM),
$$ 
Let us sum up our results.

{\bf 2.7. Theorem.} {\it The map } (2.6.2) {\it is a morphism of filtered complexes.
The induced map of $E_1$ terms for the corresponding spectral sequences is 
a pair of morphisms
$$
  \eta^i:\ DR(H_i(\fn^*, \CM(\bm)^*_N), \nabla_\KZ) \lra 
  DR(Rp_*^{N-i}\CL(\bm, N)^{\Sigma_N,-}, \nabla_\GM),
$$
$0\leq i\leq 1$. 

Here $\eta^i$ is a map from the de Rham complex of $H_i(\fn^*, \CM(\bm)^*_N)$ equipped with 
the KZ connection to the de Rham complex of $Rp_*^{N-i}\CL(\bm, N)^{\Sigma_N,-}$ equipped with 
the Gauss-Manin connection, or, which is the same, a morphism of lisse $\CD$-modules over $U_n$:
$$
  \eta^i:\ (H_i(\fn^*, \CM(\bm)^*_N), \nabla_\KZ) \lra 
  (Rp_*^{N-i}\CL(\bm, N)^{\Sigma_N,-}, \nabla_\GM).
$$ 
These maps are isomorphisms for generic $\kappa$.}

{\bf 2.8. Corollary: integral solutions for higher KZ.} Let us return to the notations of 2.1. 

Consider the complex $C^{\bullet}(\fn, M(\fm)^c)_N$, see (2.1.3)
$$
0 \lra M(\fm)_N \lra \fn\otimes  M(\fm)_{N-1} \lra 0
$$
which we denote here for brevity
$$
C^\bullet:\ 0 \lra C^0 \overset{d}\lra C^1 \lra 0,
$$
and the dual complex
$$
C^{*\bullet}:\ 0 \lra C^{1*} \overset{d^*}\lra C^{0*}\lra 0
$$
In this subsection we consider the analytic version of our varieties and $\CD$-modules.

For any $\bz = (z_1, \ldots, z_n)\in U_n$ we denote by $F(\bz)$ the fiber
$$
F(\bz) := p^{-1}(\bz) = \{(t_1, \ldots, t_N)\in \BC^N|\ t_i \neq t_j;t_i \neq z_a\}
\subset \BC^N.
$$
We will deal with the analytic Coulomb $\CD$-module $\CL^{\an}(\fm)$ over $U_{n,N}$. Consider its  de Rham complex
$$
\Omega^{\an\bullet}_{\fm} := DR(\CL^{\an}(\fm)).
$$
For each $\bz\in U_{n}$ let $\Omega_\fm^\bullet(\bz)$ denote the restriction of 
$\Omega^{\an\bullet}_{\fm}$ to the fiber $F(\bz)$;   
inside it we have the skew-symmetric part 
$$
\Omega^{\bullet}_{\fm}(\bz)^{\Sigma_N,-}\subset \Omega_\fm^\bullet(\bz)
$$
Next, inside $\Omega^{\bullet}_{\fm}(\bz)^{\Sigma_N,-}$ consider the finite-dimensional Aomoto subcomplex of differential forms with logarithmic singularities 
along all hyperplanes $t_i = t_j$ and $t_i = z_a$; let us denote this subcomplex
$$
A^\bullet(\bz):\  0 \lra A^{N-1}(\bz) \overset{d_A(\bz)}\lra A^N(\bz)\lra 0, 
$$
the differential $d_A(\bz)$ being the multiplication by the  one-form
$$
\frac{1}{\kappa}\omega_\fm(\bz) = \frac{1}{\kappa}\biggl(\sum_{1\leq i<j\leq N}2\,\dl(t_i-t_j)  
- \sum_{i=1}^N\sum_{s=1}^n m_s\,\dl(t_i-z_s)\biggr)\in \Omega^1(F(\bz))
%\eqno{(2.3.1)}
$$ 
cf. (2.3.1). 
This subcomplex will have only two nontrivial components living in degrees $N-1$ and $N$.

We denote by 
$$
W^{i}(\bz) := H^i(A^{\bullet}(\bz)),\ i = N - 1, N,
$$
its cohomology.

\

{\it Global maps $\eta^i$}

\

The space 
$$
C^0 = M(\fm)_N
$$ 
admits a base $\{f^av,\ |a| = N\}$; let us denote by $\{f^{a\vee}\}$ the dual base of 
$C^{0*}$. 

Similarly, the space 
$$
C^1 = \fn\otimes M(\fm)_{N-1}
$$ 
admits a base $\{f\otimes f^bv,\ |b| = N - 1\}$; let us denote $\{f^\vee\otimes f^{b\vee}\}$ the dual base of $C^{1*}$. 

Define two maps
$$
\eta^i:\ C^{i*} \lra \Omega^{N-i}(U_{n,N}),\ i = 0, 1,
$$
by
$$
\eta^0(f^{a\vee}) = w_a,
$$
and 
$$
\eta^1(f^\vee\otimes f^{b\vee}) = w_b.
$$
Denote
$$
\CA^i : = \eta^i(C^{i*})\subset \Omega^{N-i}(U_{n,N}).
$$
Let $\bz\in U_n$. 
The restriction to the fiber $F(\bz)$ induces maps
$$
\CA^{N-i} \lra A^{N-i}(\bz);
$$
composing them with the maps $\eta_i$ we get maps 
$$
\eta^i(\bz):\ C^{i*} \lra A^{N-i}(\bz)
$$
According to  [SV] these maps are isomorphisms; moreover, they induce 
an isomorphism of complexes
$$
\eta^\bullet(\bz): C^{*\bullet}\iso A^\bullet(\bz)[N]
$$
where on the left we have the Chevalley differential whereas on the right 
we have the twisted de Rham differential in the de Rham complex of the fiber.

\

{\it Chains of the Betti realization}

\

For each $\bz$ let $\CL_\fm(\bz)$ denote the restriction of $\CL(\fm)$ to $F(\bz)$; 
let 
$$
\CL_\fm(\bz)^\vee = \CH om(\CL_\fm(\bz),\CO^\an_{F(\bz)})
$$
be the dual $D$-module. Let 
$$
\CL_\fm(\bz)^{\vee\text{hor}}\subset \CL_\fm(\bz)^\vee
$$
be the subsheaf of horizontal sections; it isa locally constant 
sheaf over $F(z)$.

Let $C_\bullet(F(\bz), \CL_\fm(\bz)^{\vee\text{hor}})$ 
denote the complex 
$$
0 \lra 
C_{2N}(F(\bz), \CL_\fm(\bz)^{\vee\text{hor}}) \lra \ldots \lra 
C_{0}(F(\bz), \CL_\fm(\bz)^{\vee\text{hor}}) \lra 0
$$
of finite 
singular chains with coefficients in $\CL_\fm(\bz)^{\vee\text{hor}}$. We will be dealing with 
subspaces of $i$-cycles
$$
Z_{i}(F(\bz), \CL_\fm(\bz)^{\vee\text{hor}})\subset C_{i}(F(\bz), \CL_\fm(\bz)^{\vee\text{hor}})
$$
and with  homology spaces $H_{i}(F(\bz), \CL_\fm(\bz)^{\vee\text{hor}})$. 

{\bf 2.8.1. GM connection: Betty realization.} When $\bz$ varies, these complexes form a complex of (infinite dimensional) vector bundles over $U_n$, denoted by 
$\CC_{\bullet}(F(\bz), \CL_\fm(\bz)^{\vee\text{hor}})$. Each term 
$\CC_{i}(F(\bz), \CL_\fm(\bz)^{\vee\text{hor}})$ carries a flat connection. 

Indeed, 
given $\bz_0$ and a finite singular chain 
$$
\gamma(\bz_0)\in C_{i}(F(\bz_0), \CL_\fm(\bz_0)^{\vee\text{hor}}),
$$
when can move $\bz$ in a small neighbourhood $V\ni \bz_0$ such that nothing changes topologically; 
this provides a parallel transport of $\gamma(\bz_0)$ over $V$, i.e. a flat family 
of chains 
$$
\{\gamma(\bz)\in C_{i}(F(\bz), \CL_\fm(\bz)^{\vee\text{hor}})\}_{\bz\in V}.
\eqno{(2.8.1.1)}
$$
These connections are obviously compatible with boundary, i.e. we get a flat connection on 
the complex $\CC_{\bullet}(F(\bz), \CL_\fm(\bz)^{\vee\text{hor}})$. This is the Betty incarnation of the derived GM connection. 

It induces flat connections on the bundles of cycles 
$\CZ_{i}(F(\bz), \CL_\fm(\bz)^{\vee\text{hor}})$ and on the homology 
$\CH_{i}(F(\bz), \CL_\fm(\bz)^{\vee\text{hor}})$. 

$\square$ 

We can integrate $i$-forms against $i$-chains, i.e. we have pairings
$$
\int:\ C_i(F(\bz), \CL_\fm(\bz)^{\vee\text{hor}})\otimes \Omega^i_\fm(\bz) \lra \BC.
$$
Let  
$$
\{\gamma_i(\bz) \in Z_i(F(\bz), \CL_\fm(\bz)^{\vee\text{hor}})\}_{\bz\in V}
\eqno{(2.8.1)}
$$
be a flat family of cycles over a small open $V\subset U_n$
whose classes in $H_i(F(\bz), \CL_\fm(\bz)^{\vee\text{hor}})$ form a flat 
section of the GM connection.

Let $i = N$. For each $x\in C^{0*}$ and $\bz\in V$ we get a number
$$
\int_{\gamma_N(\bz)} \eta^0(\bz)(x)\in \BC;
$$
it is linear with respect to $x$, so we've got an element 
$$
\int_{\gamma_N(\bz)} \eta^0(\bz)(\bullet)\in (C^{0*})^* = C^0 = M(\fm)_N.
\eqno{(2.8.2)}
$$
Note that if $x = d^*y$ then $\eta^0(\bz)(x)$ is a coboundary in $\Omega_\fm^\an(\bz)$, 
so the integral is zero since $\gamma_N(\bz)$ is a cycle. This means that (2.8.2) 
belongs to the subspace of "singular vectors"
$$
(\Coker(d^*))^* = \Ker d = \Ker(e: M(\fm)_N\lra M(\fm)_{N-1}).
$$

Similarly if $i = N-1$ then for each $x\in C^{1*}$ and $\bz\in V$ we get a number
$$
\int_{\gamma_{N-1}(\bz)} \eta^1(\bz)(x)\in \BC;
$$
which is linear with respect to $x$, so we've got an element 
$$
\int_{\gamma_{N-1}(\bz)} \eta^1(\bz)(\bullet)\in (C^{1*})^* = C^1 = M(\fm)_{N-1}.
%\eqno{(2.8.3)}
$$
Its image in 
$$
(\Ker(d^*))^* = C^1/dC^0 = \Coker(e: M(\fm)_N\lra M(\fm)_{N-1})
$$
depends only on the homology class
$$
\overline{\gamma_{N-1}(\bz)} \in Z_i(F(\bz), \CL_\fm(\bz)^{\vee\text{hor}}).
$$
 
%NOTATION FOR THESE TWO SPACES ?

{\bf 2.8.1. Theorem.} (a) {\it 
For any local flat family of $N$-cycles 
$$
\{\gamma_{N}(\bz)\in  Z_N(F(\bz), \CL_\fm(\bz)^{\vee\text{hor}})\}_{z\in V},\ V\subset U_n, 
$$
the linear map $\int_{\gamma(\bz)}$ defines a solution of the KZ equations with values in 
$Ker (d)$, i.e. in the weight component of 
$$
\Ker e:\ M(\fm) \lra M(\fm)
$$ 
of weight $\sum_{s=1}^n m_i -2N$. For generic $\ka$ any solution of the KZ equations in this
space is given by a suitable family $\gamma_{N}(\bz)$.}  

(b) {\it 
For any local flat family of $(N-1)$-cycles 
$$
\{\gamma_{N-1}(\bz)\in Z_{N-1}(F(\bz), \CL_\fm(\bz)^{\vee\text{hor}})\}_{z\in V},\ V\subset U_n,
$$
the linear map $\int_{\gamma(\bz)}$ defines a solution of the KZ equations with values in 
$\Coker(d)$, i.e. in the weight component of 
$M(\fm)/eM(\fm)$ of weight $\sum_{s=1}^n m_i -2(N-1)$. For generic $\ka$ any solution of the KZ equations in this
space is given by a suitable family $\gamma_{N-1}(\bz)$.} 

{\bf Proof.} Part (a) is proved in [SV], whereas part (b) is new and follows from Theorem 2.7 $\square$

\

See [CV], where the dimensions of the  spaces $\Ker d$ and $\Coker d$  are
calculated for nonnegative integers $m_1,\ldots, m_n$.

\

{\bf 2.9. Exotic (dual) KZ equations.} 
%We put ourselves into the situation of [SV2], 1.1, with 

Let $N = 1, n = 2$. Let us look up more attentively at the KZ - Coulomb part of our cocycle.

% ADJUST NOTATION

So we have a $2$-dimensional subspace 
$$
M_1\subset M(m_1)\otimes M(m_2)
$$
with a base $\{ fv_1\otimes v_2,\ v_1\otimes fv_2\}$ whose elements we will write as columns.
  
The Casimir $\Omega$ acts on this subspace by the matrix
$$
\Omega =  \left(\begin{matrix} (m_1 - 2)m_2/2 & m_2\\
m_1 & m_1(m_2 - 2)/2
\end{matrix}\right)
$$ 
Consider a double Coulomb - KZ complex $\Omega^{\bullet\bullet}(M)$: as a graded 
space 
$$
\Omega^{\bullet\bullet}(M_1) := \Omega^{\bullet\bullet}(U_{2,1})\otimes M_1
$$ 
where $\Omega^{ij}$ are differential forms in $z, t$, of degree $i$ (resp. $j$)  with respect to  $z$ 
(resp. to $t$). 

The first (horizontal) differential is a KZ connection 
$$
d' = \nabla_{KZ} = d_z - \frac{1}{\kappa}\frac{\Omega(dz_1 - dz_2)}{z_1 - z_2} = 
d_z + A_1dz_1 + A_2dz_2
$$
where $d_z$ means de Rham with respect to $z$, whereas the second (vertical) differential 
$$
d'' = d_t
$$
(de Rham with respect to $t$)

The identity $\nabla_{KZ}^2 = 0$ means that the KZ connection is integrable.

In coordinates:
$$
\dpar_{z_2}A_1 - \dpar_{z_1}A_2 - [A_1, A_2] = 0,
$$
in our case $[A_1, A_2] = 0$.

Now we will descibe the relevant part of the cocycle $\CI$ from Theorem 2.5. 

Consider a form
$$
\omega^{01} = I = \left(\begin{matrix}
(t - z_1)^{-1}\Phi dt \\ (t - z_2)^{-1}\Phi dt
\end{matrix}\right) = \left(\begin{matrix}
I_1 \\ I_2
\end{matrix}\right) \in \Omega^{01}(M_1) 
$$

{\bf 2.9.1.  Claim.} {\it We have 
$$
d''\omega^{01} = 0,
$$
(obvious), and
$$
d'\omega^{01} = d''\omega^{10}
\eqno{(2.9.1)} 
$$
where 
$$
\omega^{10} = J_1 dz_1 + J_2 dz_2 \in \Omega^{10}(M_1),
$$
with
$$
J_1 = \left(\begin{matrix} - (t - z_1)^{-1}\Phi \\ 0
\end{matrix}\right),\ J_2 = \left(\begin{matrix} 0 \\ - (t - z_2)^{-1}\Phi 
\end{matrix}\right). 
$$}

{\bf 2.9.2. Claim.} {\it We have
$$
d'\omega^{10} = 0,
$$
in coordinates
$$
- \dpar_{z_2}J_1 +  \dpar_{z_1}J_2 - \frac{1}{\kappa}\frac{\Omega J_2}{z_1 - z_2} 
+ \frac{1}{\kappa}\frac{\Omega J_1}{z_1 - z_2} = 0.
\eqno{(2.9.2)}
$$}

%CHECK THE SIGN

The last differential equation is called the {\it dual KZ equation}: it is a system 
of two linear differential equations on two functions (nonzero coordinates of vectors 
$J_1, J_2$). 

The equation does not 
depend on $t$, whereas our vectors $J_1, J_2$ do. For all $t$ the couple
$$
(J_1(t, z), J_2(t, z))
$$
is a solution of (2.8.2). 

\
  
  \

  \
  
  \centerline{\bf \S 3. Kac-Moody case}
  
  \

  {\bf 3.1. Kac-Moody algebras without Serre relations.} We start with the data from [SV], 6.1. Let $\fh$ be a finite-dimensional 
  vector space equipped with a non-degenerate symmetric bilinear form $( , )$. 
  
  We fix a finite set of non-zero covectors $\{\alpha_1, \ldots, \alpha_r\}\subset \fh^*$ 
  whose elements are called simple roots; let $B = (b_{ij})$ where $b_{ij} = 
  (\alpha_i, \alpha_j)$ (this is "the symmetrized Cartan matrix"). 
  
  We denote by 
  $$
  h_i = b(\alpha_i)
  $$
  where $b: \fh^*\iso \fh$ is the isomorphism induced by $( , )$.
  
  We define $\fg = \fg(B)$ as a Lie algebra with generators $e_i, f_i, 1\leq i\leq r$, and $\fh$ and relations
  $$
  [e_i, f_j] = \delta_{ij}h_i,
  $$
  $$
  [h, e_i] = \alpha_i(h)e_i,\ [h, f_i] = - \alpha_i(h)f_i,
  $$
  $$
  [h, h'] = 0,\ h, h'\in \fh
  $$
  We denote by $\fn = \fn_-\subset \fg$ (resp. by $\fn_+$) the Lie subalgebra 
  generated by all elements $f_i$ (resp. $e_i$); it is a free Lie algebra with these generators.
  
  We have the triangular decomposition
  $$
  \fg = \fn_-\oplus \fh \oplus \fn_+
  $$
  
  % RELATION TO KM WITHOUT SERRE RELATIONS
  
  \
  
  {\it Root lattice} 
  
  \
  
  Let 
  $$
  \Lambda = \sum_i \BZ\alpha_i\subset \fh^*
  $$ 
  denote the abelian subgroup generated by $\alpha_i$. 
  
  We will use the notations for "positive" and "negative"\ submonoids:
$$
\Lambda_{\geq 0} := \sum_{i=1}^r \BZ_{\geq 0} \alpha_i \subset \Lambda,\ 
\Lambda_{> 0} := \Lambda_{\geq 0}\setminus \{\bze\};
$$
$$ 
\Lambda_{\leq 0} := - \Lambda_{\leq 0}, \Lambda_{< 0} := \Lambda_{\leq 0}\setminus \{\bze\} 
$$

\

{\it Principal gradation}

\

Our algebra $\fg$ is $\Lambda$-graded:
$$
\fg = \oplus_{\lambda\in \Lambda}\ \fg_\lambda
$$
where
$$
\fg_\lambda = \{x\in \fg|\ [h, x] = \lambda(h)x\ \text{for all}\ h\in \fh\}
$$
with
$$
\fh = \fg_\bze,\ 
$$
$$ 
\fn := \fn_- = \oplus_{\lambda\in \Lambda_{<0}}\ \fg_\lambda = \oplus_{\lambda\in \Lambda_{<0}}\ \fn_\lambda
$$
$$
\fn_+ =\oplus_{\lambda\in \Lambda_{>0}}\ \fg_\lambda = 
\oplus_{\lambda\in \Lambda_{>0}}\ \fn_\lambda
$$

\

{\it Verma modules}

\

For $\mu\in \fh^*$ $M(\mu)$ will denote a $\fg$-module with one generator $v = v_\mu$ and relations
$$
h v_\mu = \mu(h)e_\mu,\ e_iv_\mu = 0. 
$$
It is $(\mu + \Lambda_{\leq 0})$-graded:
$$
M(\mu) = \oplus_{\lambda\in \mu +  \Lambda_{\leq 0}} M(\mu)_\lambda
$$
where
$$
M(\mu)_\lambda = \{x\in M(\mu)|\ hx = \lambda(h)x\}
$$
 
A map 
$$
U\fn \lra M(\mu),\ x\mapsto xv_\mu
$$
is an isomorphism of vector spaces.

\

{\it Notation: duals for $\Lambda$-graded spaces}

\

In the sequel we will be dealing with various $\Lambda$-graded spaces 
$V = \oplus_{\lambda\in \Lambda} V_\lambda$ with finite dimentional components $V_\lambda$. 
In that case $V_\lambda^*$ will denote the restricted dual:
$$
V^* = \oplus_{\lambda\in \Lambda} V^*_\lambda.
$$ 

\

{\it Double}

\

The Borel Lie subalgebra 
$$
\fb := \fn\oplus \fh\subset \fg
$$ 
carries a structure of a {\it Lie bialgebra} (see [D]) described in 
[SV], 6.14.1. This means in particular that we have a cobracket map 
$$
\fb \lra \fb\wedge \fb
$$
which gives, after the passage to duals, a Lie algebra structure on the space $\fb^*$.  The projection 
$\fb \lra \fn$ induces an embedding $\fn^*\hra \fb^*$, and the 
subspace $\fn^*$ is a Lie subalgebra of $\fb^*$. 

This allows one to define its {\it Drinfeld double} $\tfg = D(\fb)$; it is a Lie algebra which as a vector space is
$$
\tfg = \fb \oplus \fb^* = \fn\oplus \fh \oplus \fh^*\oplus\fn^*. 
$$

If $M$ is  a Verma module,  one introduces  a structure of a $\fb^*$-module 
on $M^*$ which, together with an obvious structure of a  a $\fb$-module gives rise to a 
$D(\fb)$-module structure on $M^*$, see [SV], 6.16. 

\

{\bf 3.2. Chevalley complexes.} We fix $n\geq 1$ and an $n$-tuple of weights 
$$
\mu = (\mu_1, \ldots, \mu_n)\in \fh^{*n}
$$
Consider
$$
M(\mu) = M(\mu_1)\otimes\ldots\otimes M(\mu_n)
$$
The $\Lambda$-gradations on each $M(\mu_i)$ gives rise to a $\Lambda$-gradation on their 
tensor product $M(\mu)$.

Each $M(\mu_i)^*$ is a $\tfg$-module, whence the tensor product
$$
M(\mu)^* = M(\mu_1)^*\otimes\ldots\otimes M(\mu_n)^*
$$
is a $\tfg$-module as well. In particular due to the inclusions of Lie algebras 
$$
\fn^*\subset \fb^* \subset D(\fb^*) = \tfg
$$
$M(\mu)^*$ is a $\fn^*$-module. 

We will be interested in Chevalley homology complexes:
$$
C_\bullet(\fn^*, M(\mu)^*):\ \ldots \lra \Lambda^2\fn^*\otimes M(\mu)^*  \lra 
\fn^*\otimes M(\mu)^* \lra M(\mu)^* \lra 0
\eqno{(3.2.1)} 
$$
They are analogues of (2.1.1).

They carry a $\Lambda$-grading induced by gradings on $\fn$ and $M(\mu)$:
$$
C_\bullet(\fn^*, M(\mu)^*) = \oplus_{\lambda\in \Lambda} C_\bullet(\fn^*, M(\mu)^*)_\lambda
$$
where we denote by $C_\bullet(\fn^*, M(\mu)^*)_\lambda$ the subcomplex of weight 
$|\mu| + \lambda$, $|\mu| := \sum_{a=1}^n \mu_a$. 

\

{\bf 3.3. The Casimir element and KZ equation.}

We have an invariant Casimir element
$$
\Omega \in \tfg \hat{\otimes} \tfg^*
\eqno{(3.3.1)}
$$
defined in [SV], 7.2. Namely, 
$$
\Omega := \sum_{\lambda\in \Lambda_{<0}}\Omega_{\lambda} + \Omega_0 + \sum_{\lambda\in \Lambda_{>0}}\Omega_{\lambda}\in 
$$
$$
\in\fn\otimes\fn^*\oplus\fh\otimes\fh^*\oplus\fh^*\otimes\fh\oplus \fn^*\otimes\fn
$$
where 
$$
\Omega_0 = \frac{1}{2}(\Omega_\fh + \Omega_{\fh^*})
$$
and 
$\Omega_{\lambda}\in \fn_\lambda\otimes\fn_\lambda^*$ for $\lambda < 0$ (resp. $\in \fn^*_\lambda\otimes\fn_\lambda$ for $\lambda > 0$) are canonical elements.

Recall the space $U_n$ from 1.1. 

Let $\CC_\bullet(\fn^*, M(\mu)^*) = C_\bullet(\fn^*, \CM(\mu)^*)$ be the trivial 
vector bundle over $U_n$ with a fiber $C_\bullet(\fn^*, M(\mu)^*)$; it is a $\Lambda$-graded complex of vector bundles:
$$
C_\bullet(\fn^*, \CM(\mu)^*) = \oplus_{\lambda\in \Lambda} C_\bullet(\fn^*, \CM(\mu)^*)_\lambda
$$
For
$$
\lambda = - \sum_{i=1}^r k_i\alpha_i \in \Lambda_{\leq 0}
$$
with $\sum_{i=1}^r k_i = N$ the complex  $\CC_\bullet(\fn^*, M(\mu)^*)_{\lambda}$
leaves in degrees $[-N, 0]$.

The invariant Casimir element allows one to define the KZ connection on each 
$\CC_\bullet(\fn^*, M(\mu)^*)_\lambda$, see 1.2.

\

{\bf 3.4. A Coulomb $\CD$-module and its de Rham complex.}
Pick 
$$
\lambda = - \sum_{i=1}^r k_i\alpha_i \in \Lambda_{\leq 0};
$$
let $N = \sum k_i$. 

Consider the space $\BC^{n,N} = \BC^{n+N}$ with coordinates $z_1, \ldots, z_n, t_1, \ldots, t_N$ and a subspace 
$$
U_{n, N} = \{ (\bz, \bt)\in \BC^{n,N}|\ z_i \neq z_j, t_i \neq t_j, z_i\neq t_j\}
$$
We have a projection
$$
p = p_{n,N}:\ U_{n, N} \lra U_n
$$
We shall use a notation $[k] = \{1, \ldots, k\}$.

Pick a map of sets  
$$
\pi: [N] \lra [r]
$$
such that 
$$
|\pi^{-1}(i)| = k_i,\ i\in [r].
$$
We will denote by 
$$
\Sigma_\pi \isom \Sigma_{k_1} \times\ldots\times \Sigma_{k_r}
$$ 
a subgroup of the symmetric group respecting all the fibers 
$\pi^{-1}(i)$.

{\it A Coulomb $\CD$-module $\CL(\mu, \lambda)$}

By definition $\CL(\mu, \lambda)$ is a $\CD_{U_{n,N}}$-module which is  $\CO_{U_{n,N}}$ equipped with an integrable connection
$$
\nabla_{\mu,\lambda} = d_{DR} + \frac{1}{\kappa}\omega_{\mu,\lambda}
$$
where $\omega_{\mu,\lambda}$ is 
 a closed differential $1$-form
$$
\omega_{\mu,\lambda} = \sum_{1\leq i < j\leq n}(\mu_{i}, \mu_{j})\frac{dz_i - dz_j}{z_i - z_j} - \sum_{i\in [n], k\in [N]}(\mu_{i}, \alpha_{\pi(k)})\frac{dz_i - dt_k}{z_i - t_k} + 
$$
$$
+ \sum_{1\leq k < l \leq N}(\alpha_{\pi(k)}, \alpha_{\pi(l)})\frac{dt_k - dt_l}{t_k - t_l}
\eqno{(3.4.1)}
$$

It gives rise to the de Rham complex
$$
\Omega^\bullet_{\mu,\lambda}(U_{n,N}) := DR(\CL(\mu,\lambda))(U_{n,N})  = (\Omega^\bullet(U_{n,N}), \nabla_{\mu,\lambda})
$$

We will be interested in the subcomplex of $\Sigma_\pi$-skew-invariants
$$
\Omega^\bullet_{\mu,\lambda}(U_{n,N})^{\Sigma_\pi,-}\subset \Omega^\bullet_{\mu,\lambda}(U_{n,N})
$$

\

{\bf 3.5. Relative de Rham complexes and derived Gauss-Manin.} 

(a) The de Rham complex $\Omega^\bullet(U_{n,N})$ is the total complex of 
a bicomplex
$$
\Omega^\bullet(U_{n,N}) = \Tot\Omega^{\bullet\bullet}(U_{n,N})
$$
where $\Omega^{pq}(U_{n,N})$ is the space of forms containing $p$ differentials $dt_i$ and 
$q$ differentials $dz_m$, the full de Rham differential being the sum
$$
d_{DR} = d_z + d_t.
$$
The {\it relative} de Rham complex is by definition
$$
\Omega^\bullet(U_{n,N}/U_n) = (\Omega^{0\bullet}(U_{n,N}), d_t); 
$$
one has a projection
$$
p:\  \Omega^{\bullet}(U_{n,N}) \lra \Omega^\bullet(U_{n,N}/U_n)
$$

\

(b) {\it Coulomb twisting} 

Similarly the form 
$$
\omega_{\mu,\lambda} = \omega_{\mu,\lambda,z} + \omega_{\mu,\lambda,t} 
$$
with
$$
\omega_{\mu,\lambda,t} = \sum_{i\in [n], k\in [N]}(\mu_{i}, \alpha_{\pi(k)})\frac{dt_k}{z_i - t_k} + 
\sum_{1\leq k < l \leq N}(\alpha_{\pi(k)}, \alpha_{\pi(l)})\frac{dt_k - dt_l}{t_k - t_l}
$$
which is $d_t$-closed.

We define
$$
\Omega^\bullet_{\mu,\lambda}(U_{n,N}/U_n) := (\Omega^\bullet(U_{n,N}/U_n), d_t + \frac{1}{\kappa}\omega_{\mu,\lambda,t})
$$
% which is nothing else than $Rp_*\CL(\mu,\lambda)$

We have an epimorphism of complexes
$$
p:\ \Omega^\bullet_{\mu,\lambda}(U_{n,N}) \lra \Omega^\bullet_{\mu,\lambda}(U_{n,N}/U_n)
\eqno{(3.5.1)}
$$

(c) {\it Derived Gauss - Manin connection.}

The complex $\Omega^\bullet_\mu(U_{n,N})$ is the total complex of a double complex
$$
\Omega^{\bullet\bullet}_{\mu,\lambda}(U_{n,N}) := (\Omega^{\bullet\bullet}(U_{n,N}), \nabla_{\mu,\lambda,z} + \nabla_{\mu,\lambda,t})
$$
where
$$
\nabla_{\mu,\lambda,z} = d_z + \frac{1}{\kappa}\omega_{\mu,\lambda,z},\ 
\nabla_{\mu,\lambda,t} = d_z + \frac{1}{\kappa}\omega_{\mu,\lambda,t}.
$$
We shall write the differential $\nabla_{\mu,z}$ (resp. $\nabla_{\mu,t}$) horizontally 
(resp. vertically).

The map $p$ (3.5.1) is nothing but the projection to the utmost left vertical component.

We can identify $\Omega^{\bullet\bullet}_{\mu,\lambda}(U_{n,N})$ with the de Rham complex of the connection $\nabla_{\mu,\lambda,z}$ on 
the complex $\Omega^\bullet_{\mu,\lambda}(U_{n,N}/U_n)$:
$$
\nabla_{\mu,\lambda,z}:\ \Omega^\bullet_{\mu,\lambda}(U_{n,N}/U_n) = \Omega^{0\bullet}_{\mu,\lambda}(U_{n,N})\lra \Omega^\bullet_{\mu,\lambda}(U_{n,N}/U_n)\otimes 
\Omega^1(U_n) = \Omega^{1\bullet}_{\mu,\lambda}(U_{n,N})
$$
This is the {\it derived GM connection} on the complex $Rp_*\CL(\mu,\lambda)(U_n)$.

% EXPLAIN NOTATION 

{\bf 3.7. A map $\eta$ and its lifting $\teta$. }

In [SV] a map of complexes
$$
\eta: C_\bullet(\fn^*, \CM(\mu)^*)_\lambda \lra \Omega^\bullet_{\mu,\lambda}(U_{n,N}/U_n)^{\Sigma_\pi,-}[N]
\eqno{(3.7.1)}
$$
has been defined, see {\it op. cit.} (7.2.4). 

Here we consider both complexes appearing in (3.7.1) as 
cohomological complexes concentrated in degrees $[-N, 0]$.

For each $\bz\in U_n$ consider a fiber
$$
U_\bz = U_{n,N;\bz} := p^{-1}(\bz)
$$
We can compose $\eta$ with the restriction map
$$
r_\bz: \Omega^\bullet_{\mu,\lambda}(U_{n,N}/U_n)\lra \Omega^\bullet_{\mu,\lambda}(U_{n,N;\bz})
$$
to get 
$$
\eta_\bz = r_\bz\circ\eta: C_\bullet(\fn_-, \CM(\mu)^*)_\lambda \lra 
\Omega^\bullet_{\mu,\lambda}(U_{n,N;\bz})^{\Sigma_\pi,-}[N]
$$

A remarkable feature of the mappings $\eta_\bz$ is the following:

{\it for generic values of $\kappa$ the maps $\eta_\bz$ are quasi-isomorphisms for all 
$\bz\in U_n$.}

Here "generic"\ means $\kappa\in\BC\setminus$ (an explicitly given discrete countable subset).

\

{\it Main result}

\

We start with 
is a definition of  a map of {\it graded $\CO(U_n)$-modules} 
$$
\teta: C_\bullet(\fn^*, \CM(\mu)^*)_\lambda \lra 
\Omega^\bullet_{\mu,\lambda}(U_{n,N})^{\Sigma_\pi,-}[N]
\eqno{(3.7.2)}
$$
which lifts $\eta$, i.e. such that
$$
\eta = p\circ\teta.
$$
Here is a picture:
$$
\begin{matrix} & & \Omega^\bullet_{\mu,\lambda}(U_{n,N})^{\Sigma_\pi,-}[N]\\
& \teta\nearrow & \downarrow p\\
C_\bullet(\fn^*, \CM(\mu)^*)_\lambda & \overset\eta\lra & 
\Omega^\bullet_{\mu,\lambda}(U_{n,N}/U_n)^{\Sigma_\pi,-}[N]
\end{matrix}
$$
The definition of $\teta$ is a modification of that of $\eta$. Namely, 
for a monomial
$$
x\in C_i(\fn^*, \CM(\mu)^*)_\lambda
$$
the corresponding differential form 
$$
\eta(x)\in \Omega^{N-i}_{\mu,\lambda}(U_{n,N}/U_n)
$$
contains fractions of the form $dt_i/(t_i - z_p)$. To obtain $\teta(x)$ we replace 
all these fractions by $d\ln(t_i - z_p)$. That's it.

Confer the definition of forms $u_a, u_b$ in 2.2 for $\fg = \fsl_2$.

The map $\teta$ induces a map of $\Omega^\bullet(U_{n,N})$-modules
$$
\teta: \Omega^\bullet(U_{n,N})\otimes_{\Omega^\bullet(U_{n})}C_\bullet(\fn^*, \CM(\mu)^*)_\lambda \lra 
\Omega^\bullet_{\mu,\lambda}(U_{n,N})^{\Sigma_\pi,-}[N]
\eqno{(3.7.3)}
$$
The space on the left is the underlying space of the  
De Rham complex for the derived KZ; it carries the KZ differential $\nabla_{KZ}$.

{\bf 3.8. Theorem.} (a) {\it The map $\teta$} (3.7.3) {\it commutes with the differentials on both sides.}

(b) {\it Both sides of} (3.7.3) {\it carry natural decreasing filtrations in $z$ direction, and the map $\teta$ respects these filtrations. } 

In other words, we've got a map of filtered complexes

\hspace{2cm} 
$$
\begin{matrix}DR(\text{derived KZ}) &\overset\teta\lra & DR(\text{derived GM})\\
\parallel & & \parallel\\
\Tot\Omega_{KZ}^\bullet(U_n, C_\bullet(\fn^*, \CM(\mu)^*)_\lambda) & 
\overset\teta\lra & \Tot(\Omega^{\bullet\bullet}_{\mu,\lambda})(U_{n,N})^{\Sigma_\pi, -}[N] 
\end{matrix}
\eqno{(3.8.1)}
$$

The proof is similar to that of Section 2.

{\bf 3.9. Corollary.} {\it The map $\teta$ induces maps of 
$\CD_{U_n}$-modules (which are isomorphisms for generic $\kappa$)
$$
\eta_i:\ (\CH_i(\fn^*, \CM(\mu)^*)_\lambda, \nabla_{KZ}) \lra ( R^{N-i}p_*\CL(\mu,\lambda)^{\Sigma_\pi, -}, \nabla_{GM})
$$
for all $0\leq i\leq N$. }

For $i = 0$ such a mapping has been constructed in [SV].

% INTEGRAL REPS OF SOLUTIONS

As in Section 2, we get from this an integral representation for the solutions.

% CHECK AND EXPLAIN NOTATION

% COROLLARY FOR IRREDUCIBLE REPS?

\

\centerline{\bf Appendix}

\

We recall here some standard constructions from homological algebra.

{\bf A.1. Bicomplexes.} A a {\it bicomplex } in an abelian category $\CC$ is a collection 
of objects $A^{\bullet\bullet} = \{ A^{pq},\ p, q\in \BZ\}$, and arrows 
$$
d_h^{pq}:\ A^{pq} \lra A^{p+1,q},\ d_h^{p,q+1}:\ A^{pq} \lra A^{p,q+1} 
$$
such that
$$
d_h^2 = d_v^2 = 0,\ d_hd_v = d_vd_h
$$
One associates to it a simple complex, to be denoted
$$
A^\bullet = \Tot A^{\bullet\bullet}
$$
with components 
$$
A^i = \oplus_{p + q = i} A^{pq}
$$
and a differential $d: A^i\lra A^{i+1}$ with components
$$
d^{pq} = d_h^{pq} + (-1)^pd_v^{pq}:\ A^{pq}\lra A^{p+1,q}\oplus A^{p,q+1}.
$$

{\bf A.2. Filtered complexes.} Let $A^\bullet$ be a simple complex.  Consider a decreasing filtration by subcomplexes on it:
$$
F^0A^\bullet = A^\bullet \supset F^1A^\bullet \supset \ldots .
$$
We associate to it a collection of complexes 
$$
E(A^\bullet, F)^i: \ 0 \lra H^i(F^0A^\bullet/F^1A^\bullet) \lra H^{i+1}(F^1A^\bullet/F^2A^\bullet) \lra \ldots,
\eqno{(A.2.1)} 
$$
$i\geq 0$, 
where a differential 
$$
H^{i+p}(F^pA^\bullet/F^{p+1}A^\bullet) \lra H^{i+p+1}(F^{p+1}A^\bullet/F^{p+2}A^\bullet)
$$
is the boundary map for the  short exact sequence
$$
0 \lra F^{p+1}A^\bullet/F^{p+2}A^\bullet \lra F^{p}A^\bullet/F^{p+2}A^\bullet \lra 
F^{p}A^\bullet/F^{p+1}A^\bullet \lra 0
$$
(This is nothing else but the $E_1$ term of the spectral sequences for 
$(A^\bullet, F^\bullet)$.) 

\

{\bf A.3. Example.} Suppose that $A^\bullet = \Tot A^{\bullet\bullet}$ 
with $A^{pq} = 0$ for $p < 0$, and a filtration 
is defined by
$$
F^iA^j = \oplus_{p\geq i, p+q = j} A^{pq}.
$$
Then a $p$-th graded piece  
$$
F^pA^\bullet/F^{p+1}A^\bullet = \{ A^{pq},\ q\in \BZ\},
$$
and the differential induced by $d$ on it coincides with the vertical differential $d_v$. 

It follows that a complex 
$E(A^\bullet, F^\bullet)^i$ is identified with 
$$
0 \lra H^i_v(A^{\bullet 0}) \lra H^i_v(A^{\bullet 1}) \lra \ldots,   
$$
with a differential induced by $d_h$.

\bigskip\bigskip

\centerline{\bf References}

\bigskip\bigskip

[CV] D.Cohen, A.Varchenko,   
Resonant local systems on complements of discriminantal arrangements and $\fsl_2$ 
representations, {\it Geom. Dedicata} {\bf 101} (2003), 217--233. 

%[DG] V.Drinfeld, D.Gaitsgory, On some finiteness questions for algebraic stacks.

[D]  V.Drinfeld, Quantum groups, Proceedings of the International Congress of Mathematicians, Vol. 1, 2 (Berkeley, Calif., 1986), 798–820, AMS, Providence, RI, 1987. 

[G] A.Grothendieck, Crystals and de Rham cohomology of schemes, Notes by I. Coates and O. Jussila, {\it Adv. Stud. Pure Math.} {\bf 3}, Dix exposés sur la cohomologie des schémas, 306–358, North-Holland, Amsterdam, 1968. 

[K] V.G.Kac, Infinite dimensional Lie algebras, Cambridge University Press, Cambridge, 1990. 

% [K] M.Kapranov, On DG-modules over the de Rham complex and the vanishing cycles functor. Algebraic geometry (Chicago, IL, 1989), 57–86,
% Lecture Notes in Math., 1479, Springer, Berlin, 1991. 

[KO] N.Katz, T.Oda, On the differentiation of de Rham cohomology classes with respect to parameters.
{\it J. Math. Kyoto Univ.} {\bf  8} (1968), 199–213. 

[SV] V.Schechtman, A.Varchenko, Arrangements of hyperplanes and Lie algebra homology, {\it 
Invent. Math.} {\bf 106} (1991), 139–194.  

%[SV2] V.Schechtman, A.Varchenko, Solutions of KZ equations modulo $p$

[S] J.-P. Serre, Lie algebras and Lie groups, 
1964 lectures given at Harvard University. Second edition. Lecture Notes in Mathematics, {\bf 1500}, Springer-Verlag, Berlin, 1992.

\newpage

%\input ???

%  WRITE HERE

\end{document}